\documentclass{article}
\usepackage[utf8]{inputenc}

\usepackage{amssymb}
\usepackage{amsmath}
\usepackage{amsfonts} 
\usepackage{amsbsy}
\usepackage{amscd}
\usepackage{amsthm}
\usepackage{geometry}                
\usepackage{graphicx}
\usepackage{amssymb}
\usepackage{epstopdf}
\usepackage{mathtools}
\usepackage{mathrsfs}
\usepackage[capitalize]{cleveref}
\usepackage{mathrsfs} 
\usepackage{subcaption}
\usepackage[framemethod=tikz]{mdframed} 
\usepackage{circuitikz}
\usepackage{overpic}
\usetikzlibrary{arrows}
\usetikzlibrary{patterns}

\newtheorem{theorem}{Theorem}[section]
\newtheorem{corollary}[theorem]{Corollary}

\newtheorem{definition}[theorem]{Definition}

\newtheorem{lemma}[theorem]{Lemma}

\newtheorem{example}[theorem]{Example}

\newtheorem*{maintheorem}{Theorem}

\newcommand{\C}{\mathbb{C}}

\newcommand{\M}{\mathcal{M}}

\newcommand{\W}{\mathbb W}

\newcommand{\Hidden}[1]{}

\newcommand{\RS}{\mathcal R_S}

\newcommand{\1}{\mathbf 1}
\def\big{\bigskip}

\title{Characterizing cospectral vertices via isospectral reduction}
\author{Mark Kempton\footnote{Deparment of Mathematics, Brigham Young University, Provo UT, USA, mkempton@mathematics.byu.edu}, John Sinkovic\footnote{Deparment of Mathematics, Brigham Young University, Provo UT, USA, sinkovic@mathematics.byu.edu}, Dallas Smith\footnote{Deparment of Mathematics, Brigham Young University, Provo UT, USA, dallas.smith@mathematics.byu.edu}, and Benjamin Webb\footnote{Deparment of Mathematics, Brigham Young University, Provo UT, USA, bwebb@mathematics.byu.edu}}

\date{}

\tikzstyle{every node}=[circle,fill=black,inner sep=1pt]

\begin{document}

\maketitle

\begin{abstract}
    Two emerging topics in graph theory are the study of cospectral vertices of a graph, and the study of isospectral reductions of graphs.  In this paper, we prove a fundamental relationship between these two areas, which is that two vertices of a graph are cospectral if and only if the isospectral reduction over these vertices has a nontrivial automorphism. It is well known that if two vertices of a graph are \emph{symmetric}, i.e. if there exists a graph automorphism permuting these two vertices, then they are cospectral. This paper extends this result showing that any two cospectral vertices are symmetric in some reduced version of the graph. We also prove that two vertices are strongly cospectral if and only if they are cospectral and the isospectral reduction over these two vertices has simple eigenvalues. We further describe how these results can be used to construct new families of graphs with cospectral vertices. 
\end{abstract}

\section{Introduction}

The basic question of spectral graph theory is: What can the eigenvalues of a matrix tell us about its underlying graph?  Fundamental to this question is the construction of cospectral graphs---that is, graphs which are not necessarily isomorphic whose adjacency matrices have the same eigenvalues.  Closely related to the notion of cospectral graphs is the idea of cospectral vertices.  We say that two vertices $a$ and $b$ of a graph $G$ are \emph{cospectral} if the graphs $G\backslash a$ and $G\backslash b$ are cospectral graphs.  Aside from their inherent combinatorial interest, cospectral vertices have been shown to have important application in the field of quantum information theory. Particularly, it has been shown that a necessary condition for the perfect transfer of a quantum state between two vertices of a network of interacting qubits is that those vertices be cospectral (see \cite{godsil_survey,kay2010perfect}).  As such, cospectral vertices have seen increased interest in recent literature on quantum computing \cite{godsil_survey,coutinho2016spectrally,godsil2017,eisenberg2018pretty,chan2019quantum}. 

Seemingly unrelated to cospectral vertices, an isospectral graph reduction is a method of reducing a graph to a graph with a smaller number of vertices while maintaining the graph's eigenvalues and eigenvectors \cite{thebook,Bunimovich_2011,PedroDuarte_2015}. As not to violate the fundamental theorem of algebra this smaller graph's edges are weighted with rational functions. Equivalently, an isospectral reduction of a matrix is a way of taking a matrix and constructing a smaller matrix whose entries are rational functions, in a way that preserves the matrix' spectral properties. Isospectral reductions have been used to improve the eigenvalue approximations of Gershgorin, Brauer, and Brualdi \cite{Gershgorin_1931,Brauer_1947,Brualdi_1982}; study the pseudo-spectra of graphs and matrices \cite{Vasquez_2015}; create stability preserving transformations of networks \cite{Bunimovich_2011,Bunimovich_2013,Reber_2019}; and study the survival probabilities in open dynamical systems \cite{bunimovich_2014}. 

Isospectral graph reductions have also been used to study ``hidden symmetries" in real and theoretical networks \cite{smith2019hidden}. A set of vertices in a simple graph $G$ is referred to as being \emph{symmetric} or \emph{automorphic} if there is a nontrivial automorphism $\phi:V\rightarrow V$ of the graph's vertices that fixes this set. More intuitively, a graph automorphism describes how pieces of a graph can be permuted in a way that preserves the graph’s overall structure. 
As this notion of symmetry can be extended to weighted graphs, a \emph{hidden} or \emph{latent symmetry} is an automorphism of some set of vertices in an isospectral reduction of a graph. Although these vertices may not be automorphic in the original unreduced graph, they share a number of properties that are known to hold for standard symmetries (for details see \cite{smith2019hidden}).

On the surface, the study of cospectral vertices and the study of isospectral reductions seem unrelated.  In this paper we prove that they are connected in a fundamental way. In particular, we prove the following theorem.

\begin{maintheorem} \textbf{(Symmetry of Cospectral Vertices)}
Let $a$ and $b$ be two vertices of a graph $G$, and let $\mathcal R = \mathcal R_{S}(G)$ be the isospectral reduction over the vertex set $S=\{a,b\}$.  Then\\
(i) $a$ and $b$ are cospectral if and only if $\mathcal R_{a,a} = \mathcal R_{b,b}$; and\\ 
(ii) if $G$ is undirected then (i) is equivalent to the statement: $a$ and $b$ are cospectral if and only if $a$ and $b$ are symmetric in the isospectral reduction $\mathcal R$.
\end{maintheorem}

In other words, two vertices are cospectral if and only if they are latently symmetric. (See Section \ref{sec:cosp} for the precise statement and proof of this result.) It is clear that if two vertices $a$ and $b$ of a graph $G$ are symmetric, then they are cospectral, since $G\backslash a$ and $G\backslash b$ are actually isomorphic (see \cite{kempton2017pretty}). One of the main contributions of this paper is to show that any two cospectral vertices are symmetric in some reduced version of the graph, i.e. cospectrality can always be thought of in terms of graph symmetries.  

In addition to cospectrality, an important property in the study of quantum state transfer is a strengthening of the cospectral condition called strong cospectrality \cite{godsil_survey}.  Two vertices $a$ and $b$ are said to be \emph{strongly cospectral} if they are cospectral and $\varphi(a) = \pm\varphi(b)$ for any eigenvector $\varphi$ of the adjacency matrix.  In addition to the above characterization of cospectrality, we also characterize strong cospectrality in terms of isospectral reductions.  Namely, we prove the following theorem (the proof is found in Section \ref{sec:strcosp}).
\begin{maintheorem} \textbf{(Isospectral Reductions and Strong Cospectrality)}
Two cospectral vertices $a$ and $b$ of a graph $G$ are strongly cospectral if and only if they are cospectral and the isospectral reduction $\mathcal{R}_S(G)$ over the vertex set $S=\{a,b\}$ has simple eigenvalues.
\end{maintheorem}
 
There is an important characterization of cospectral vertices which states that two vertices of a graph are cospectral if and only if the number of closed walks in the graph at each vertex of any given length is the same (see \cite{godsil2017}). Thus cospectrality can be understood entirely in terms of the combinatorics of the graph.  We also give an interpretation of the isospectral reduction as a generating function counting walks in a graph that do not visit a particular subset of the vertex set, except at the first and last step, or ``non-returning" closed walks.  Our results can thus be interpreted combinatorially as relating the enumeration of all closed walks at a vertex to the enumeration of non-returning closed walks at a vertex. Indeed, we find that the generating functions for all closed walks and all non-returning closed walks are related in a simple way. Section \ref{sec:walk} discusses this combinatorial interpretation.  

We note that the study of cospectral vertices has primarily been done only on undirected graphs.  Indeed, many of the characterizations of cospectrality require the matrices involved to be symmetric (see Theorem 3.1 of \cite{godsil2017}).  This is especially true in the study of strong cospectrality.  However, isospectral reductions can be done on a symmetric or asymmetric matrix, and much of the theory of isospectral reductions is actually more natural in the setting of directed graphs.  Thus, one of the contributions of this work is to extend the notion of strong cospectrality to directed graphs.  
In Section \ref{sec:ex}, we discuss how our results can be applied to constructing new families of graphs with cospectral vertices (a generally difficult problem) by ``unpacking" the isospectral reduction---that is, we look at possible ways of reverse engineering an isospectral reduction with symmetry to produce a graph with a pair of cospectral vertices.  We also discuss the notion of ``measure of latency" found in \cite{smith2019hidden} and how it relates to our work.  We suggest several possible directions of future research.

\section{Preliminaries}\label{sec:pelim}

The main mathematical objects we consider in this paper are graphs. The type of graphs we consider, for the sake of generality, are weighted graphs $G=(V,E,\omega)$ composed of a finite \emph{vertex set} $V=\{1,2,\dots,n\}$, a finite \emph{edge set} $E$, and a function $\omega:E\rightarrow\C$ used to weight the edges $E$ of the graph. In this framework, an \emph{unweighted graph} is one in which every edge is given unit weight. 

A graph can also be either \emph{directed} or \emph{undirected} depending on whether the graph's edges have a specified direction. If $G$ is directed then we let $e_{ij}\in E$ denote the \emph{directed edge} from vertex $i$ to vertex $j$ in the graph. The edge $e_{ii}\in E$ that begins and ends at vertex $i$ is referred to as a \emph{loop}. If $G$ is an undirected graph then it can be considered to be a directed graph in which every undirected edge is replaced with two directed edges oriented in opposite directions and each loop is replaced by a single directed loop.

Under these conventions the class of graphs we consider are those weighted directed graphs $G=(V,E,\omega)$ that may or may not have loops. We let $M=M(G)$ denote the graph's \emph{weighted adjacency matrix} given by
\[M_{i,j}=
\begin{cases}
\omega(i,j) &\text{  if } e_{ij}\in E\\
0 &\text{otherwise}.
\end{cases}
\]
Conversely, given a matrix $M\in\C^{n\times n}$, the \emph{graph of $M$}, denoted by $G=G(M)$, is the graph $G$ with weighted adjacency matrix $M$. 

It is worth noting that there is a one-to-one correspondence between the graphs we consider and the matrices $M\in\C^{n\times n}$. This allows us to consider graphs and their associated weighted adjacency matrices interchangeably. At times we will use notation and terminology from graph theory in reference to matrices and vice versa. For instance, by $M\backslash a$, we mean the matrix obtained by deleting the row and column indexed by $a$, so if $M=M(G)$ then $M\backslash a$ is the weighted adjacency matrix of $G\backslash a$.  Similarly, by $V(M)$ we mean the ``vertex set" of the matrix $M$, i.e., the vertices of the graph of $M$, or the index set for the rows and columns of $M$, and the characteristic polynomial of a graph is the characteristic polynomial of its associated matrix, which we denote by $p(M,\lambda) = \det(M - \lambda I)$.   

To be able to state the main results of this paper we need the notion of a graph automorphism. 

\begin{definition}\label{def:sym} \textbf{(Graph Automorphism)}
    An \emph{automorphism} $\phi$ of graph $G=(V,E,\omega)$ is a permutation of the graph's vertex set $V$ such that the graph's weighted adjacency matrix $M=M(G)$ satisfies $M_{ij} = M_{\phi(i) \phi(j)}$ for each pair of vertices $i$ and $j$. We let $Aut(G)$ denote the the graph's group of automorphisms.
\end{definition}

Importantly, a graph’s group of automorphisms characterizes the symmetries in the graph’s structure. For a graph $G=(V,E,\omega)$ with automorphism $\phi$ the relation $\sim$ on $V$ given by $u\sim v$ if and only if $v=\phi^j(u)$ for some nonnegative integer $j$ is an equivalence relation on $V$. These equivalence classes are called the \emph{orbits} of $\phi$. Here we say that two vertices are \emph{automorphic} if they belong to the same orbit under some $\phi\in Aut(G)$. 

\subsection{Cospectral and Strongly Cospectral Vertices}\label{subsec:prelim_cospectral}

Our main goal in this paper is relating graph symmetries, described by the graph's group of automorphisms, to the notion of cospectrality.

\begin{definition}\textbf{(Cospectral)}
Two vertices $a,b\in V$ of a graph $G=(V,E,\omega)$ are said to be \emph{cospectral} if the (weighted) adjacency matrices of $G\backslash a$ and $G\backslash b$ have the same characteristic polynomial.
\end{definition}

Many characterizations of cospectral vertices exist in the literature.  To state them, we need some notation.  For a symmetric matrix $M\in\mathbb{R}^{n\times n}$, we write its spectral decomposition as
\[
M = \sum_\theta \theta E_\theta
\]
where the sum is taken over all distinct eigenvalues $\theta\in\sigma(M)$. Here $E_\theta$ is the spectral idempotent for the projection onto the eigenspace corresponding to eigenvalue $\theta$ in this decomposition.  We let the matrix $W(M,t)$ denote the walk generating function defined by $M$, that is
\[
W(M,t) = \sum_{k=0}^\infty M^kt^k.
\]
Observe that the coefficient of $t^k$ in the $(a,b)$ entry of $W(M,t)$ is the (weighted) number of walks from $a$ to $b$ of length $k$ in the (weighted) graph $G=G(M)$. We point out for a vertex $a\in V$ of the graph $G=(V,E,\omega)$ that 
\begin{align}\label{eq:walk}
W\left(M,\frac1\lambda\right)_{a,a} = -\lambda\frac{p(M\backslash a,\lambda)}{p(M,\lambda)}.
\end{align}
as described in \cite{godsil2017}. With this notation in place, we can state some important known characterizations of cospectral vertices. 

\begin{lemma}[Theorem 3.1 of \cite{godsil2017}]\label{lem:cosp_char}
 Let $M$ be a symmetric matrix, and let $a,b \in V(M)$. 
 The following are equivalent:
\begin{enumerate}
\item  $a$ and $b$ are cospectral.
\item $(E_\theta)_{a,a} = (E_\theta)_{b,b}$ for all $\theta$. 
\item $(M^k)_{a,a} = (M^k)_{b,b}$  for all integers $k\geq 0$.
\item $W(M,t)_{a,a}=W(M,t)_{b,b}$. 
\end{enumerate}
\end{lemma}

We remark that in item 3 in this lemma, the entries of $M^k$ count (weighted) walks of length $k$ in $G$, so this is saying that $a$ and $b$ are cospectral if and only if the number of closed walks at $a$ and $b$ of any length are equal. 

It is interesting to ask how pairs of cospectral vertices in graphs can arise.  The most obvious way to obtain a cospectral pair is via symmetry (see \cite{kempton2017pretty,godsil2017}).  Note that if the graph $G=(V,E,\omega)$ has an automorphism that interchanges $a,b\in V$, i.e. $a$ and $b$ are automorphic, then $G\backslash a$ and $G\backslash b$ will be isomorphic as graphs. Thus $a$ and $b$ will be cospectral vertices.  However, there need not be an automorphism taking $a$ to $b$ for $G\backslash a$ and $G\backslash b$ to be isomorphic (see Figure \ref{fig:isom_no_aut}).

In addition, $G\backslash a$ and $G\backslash b$ need not be isomorphic in order for $a$ and $b$ to be cospectral, it is only necessary for their adjacency matrices to have the same spectrum.  It is known that in strongly regular graphs, any pair of vertices is cospectral \cite{godsil2017state}, and most examples of strongly regular graphs have no automorphisms and will not have $G\backslash a$ isomorphic to $G\backslash b$.  In addition, graphs with cospectral pairs can be constructed from cospectral graphs involved in Godsil-McKay switching (see \cite{godsil1982constructing}).  In addition to these standard constructions, there are a number of ad hoc examples that do not fit into any of these categories.  See, for instance, the graph in Figure \ref{fig:non-aut_cosp}.  In Section \ref{sec:unpack}, we discuss how our results might be used to construct more examples of cospectral pairs outside these known constructions.  

In addition to cospectrality, a stronger notion, called strong cospectrality, is important in the study of quantum information transfer, and construction of families of strongly cospectral vertices is also an important problem \cite{godsil2017}.  We define strong cospectrality in the following.  

\begin{figure}
    \begin{center}
\begin{tikzpicture}
\draw \foreach \x in {0,1,2,3,4}
{
(\x,0)node{}--(\x+1,0)node{}
}
(2,0)--(2.5,1)node{}--(3,0)
(4,0)--(4.5,1)node{}--(5,0)
(2,-.1)node[fill=white, below]{\small $a$}
(4,-.1)node[fill=white, below]{\small $b$};
\end{tikzpicture}
\end{center}
    \caption{A graph for which $G\backslash a$ and $G\backslash b$ are isomorphic, but no automorphism takes $a$ to $b$.}
    \label{fig:isom_no_aut}
\end{figure}
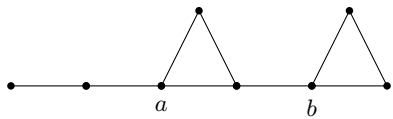

\begin{figure}
    \begin{center}
\begin{tikzpicture}
\draw \foreach \x in {0,1,2,3,4,5,6}
{
(\x,0)node{}--(\x+1,0)node{}
}
(5,0)--(5,1)node{}
(3,-.1)node[fill=white, below]{\small $a$}
(6,-.1)node[fill=white, below]{\small $b$};
\end{tikzpicture}~~~~
\begin{tikzpicture}
\draw (0,0)node{}--(-1,.6)node{}--(-1,-.6)node{}--(0,0)--(1,1)node{}--(2.1,.6)node{}--(2.1,-.6)node{}--(1,-1)node{}--(0,0)
(1,1)--(1,2)node{}
(2.1,.6)node{}--(3,0)node{}--(2.1,-.6)node{}
(-1.1,.6)node[fill=white, left]{\small $a$}
(3.1,0)node[fill=white, right]{\small $b$};
\end{tikzpicture}
\end{center}
    \caption{Two examples of graphs (both denoted $G$) with a pair of cospectral vertices $a$ and $b$, but where $G\backslash a$ and $G\backslash b$ are not isomorphic graphs, but have the same spectrum.}
    \label{fig:non-aut_cosp}
\end{figure}
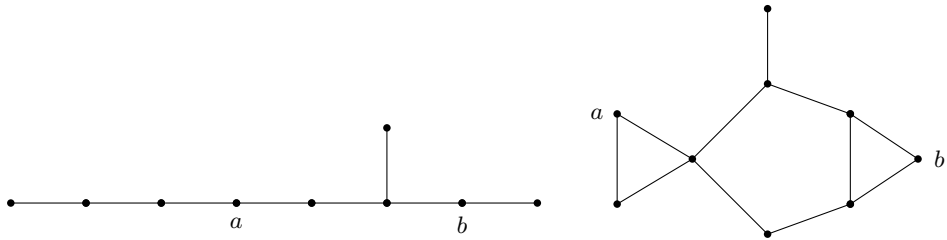

\begin{definition}\label{def:str-cosp}\textbf{(Strongly Cospectral)}
Let $e_u$ denote the characteristic vector for any vertex $u$. We say that vertices $a$ and $b$ are \emph{strongly cospectral} if
\[
E_\theta e_a = \pm E_\theta e_b
\]
for all projections $E_\theta$ in the spectral decomposition of $M$.
\end{definition}

Note that since each $E_\theta$ is symmetric, then strongly cospectral does in fact imply cospectral (see Lemma \ref{lem:cosp_char}), and is in fact a stronger notion.  Strong cospectrality is equivalent to saying that $\varphi(a)=\pm\varphi(b)$ for all eigenvectors $\varphi$ of the adjacency matrix (see \cite{kempton2017pretty}).

\subsection{Isospectral Reduction and Hidden Symmetries}\label{sec:2.2}

To formally describe the relationship between graph symmetries and cospectral vertices we need the notion of an isospectral graph reduction. An isospectral graph reduction is a method which produces a smaller graph with essentially the same set of eigenvalues as the original graph. This method for reducing a graph can be formulated both for the graph and equivalently for its associated (weighted) adjacency matrix, i.e. an isospectral graph reduction and an isospectral matrix reduction, respectively. Both types of reductions will be useful to us. For the sake of simplicity we begin by defining an isospectral matrix reduction. 

For this reduction we need to consider matrices of rational functions. The reason is that, by the Fundamental Theorem of Algebra, a matrix $M\in\mathbb{C}^{n\times n}$ has exactly $n$ eigenvalues including multiplicities. In order to reduce the size of a matrix while at the same time preserving its eigenvalues we need something that carries more information than scalars. The objects we will use to preserve this information are rational functions. The specific reasons focusing on rational functions can be found in \cite{thebook}, Chapter 1.

For a matrix $M\in\mathbb{C}^{n\times n}$ let $V=\{1,...,n\}$. If the sets $R,C\subset V$, we denote by $M_{RC}$ the $|R|\times |C|$ submatrix of $M$ with rows indexed by $R$ and columns indexed by $C$. We denote by $\bar{S}$ the complement of $S$ in $V$.  We let $\mathbb{W}^{n\times n}$ be the set of $n\times n$ matrices whose entries are rational functions $p(\lambda)/q(\lambda)\in\mathbb{W}$ with complex coefficients, where $p(\lambda)$ and $q(\lambda)\neq 0$ are polynomials with complex coefficients in the variable $\lambda$ with no common factors and $deg(p)\leq deg(q)$. The isospectral reduction of a square complex-valued matrix is defined as follows.

\begin{definition}\label{def:isored} \textbf{(Isospectral Matrix Reduction)}
The \emph{isospectral reduction} of a matrix $M\in\mathbb{C}^{n\times n}$ over the proper subset $S\subset V$ is the matrix
\begin{equation}\label{IMR}
\mathcal{R}_S(M) = M_{SS} - M_{S\bar{S}}(M_{\bar{S}\bar{S}}-\lambda I)^{-1} M_{\bar{S}S}\in\mathbb{W}^{|S|\times|S|}.
\end{equation}
where $\lambda$ is an indeterminate variable.
\end{definition}

The fact that an isospectral reduction of a matrix $M\in\mathbb{C}^{n\times n}$ results in a matrix $\mathcal{R}_S(M)\in\mathbb{W}^{|S|\times|S|}$ is shown in Theorem 1.2 of \cite{thebook}. We note that the eigenvalues of a square matrix $R=R(\lambda)\in\mathbb{W}^{n\times n}$ are defined to be solutions to its characteristic equation
\[
\det(R(\lambda)-\lambda I)=0,
\]
which is an extension of the standard definition of the eigenvalues for a matrix with complex entries. By way of notation, we let $\sigma(M)$ denote the set of eigenvalues of the matrix $M$ including multiplicities. An important aspect of an isospectral reduction is that the eigenvalues of the matrix $M$ and the eigenvalues of its isospectral reduction $\mathcal{R}_S(M)$ are essentially the same, as described by the following theorem.

\begin{theorem}[Theorem 2.2 of \cite{thebook}]\label{thm:maintheorem}\textbf{(Spectrum of Isospectral Reductions)} For $M\in\mathbb{C}^{n\times n}$ and a proper subset $S\subseteq V$, the eigenvalues of the isospectral reduction $\mathcal{R}_S(M)$ are
$$\sigma\big(\mathcal{R}_S(M)\big)=\sigma(M)-\sigma(M_{\bar{S}\bar{S}}).$$
\end{theorem}

That is, when a matrix $M$ is isospectrally reduced over a set $S$, the set of eigenvalues of the resulting matrix is the same as the set of eigenvalues of the original matrix $M$ after removing any elements which are eigenvalues of the submatrix $M_{\bar{S}\bar{S}}$. In fact, the matrix $M$ and the submatrix $M_{\bar{S}\bar{S}}$ often have no eigenvalues in common, in which case the spectrum is unchanged by the reduction, i.e. $\sigma(\mathcal{R}_S(M))=\sigma(M)$.

A fact that will be important to us is that an isospectral matrix reduction preserves the symmetry of a symmetric matrix (see Theorem 2.2 and \cite{thebook}).
 
\begin{lemma}\label{lem:symm}
If $M\in\mathbb{C}^{n\times n}$ is a symmetric matrix, then $\RS(M)$ is symmetric for any $S\subset V$.
\end{lemma}
 
If the graph $G=(V,E,\omega)$ has the (weighted) adjacency matrix $M=M(G)$ then the isospectral reduction of $G$ over $S\subset V$ is defined as the graph with adjacency matrix $\mathcal{R}_S(M)$. Hence, this graph has edges that are weighted by elements of $\mathbb{W}$.  We call this graph the \emph{isospectral graph reduction} of $G$ over $S$ and denote it by $\mathcal{R}_S(G)$.  

We note that although an isospectral graph reduction is defined in terms of an isospectral matrix reduction, it can be obtained entirely via graph theoretic operations.  This perspective will be useful to us as well, so we present it in the following section.  
\subsubsection{Isosepctral Graph Reductions}

An isospectral graph reduction can be described in terms of collapsing certain graph paths and cycles to single edges and loops, respectively. Recall that a \textit{path} $P$ in the graph $G=(V,E,\omega)$ is an ordered sequence of distinct vertices $P=1,\dots,m\in V=\{1,2,\dots,n\}$ such that $e_{i,i+1}\in E$ for $1\leq i\leq m-1$. We call the vertices $2,\dots,m-1$ of $P$ the \textit{interior} vertices of $P$. If the vertices $1$ and $m$ are the same then $P$ is a \textit{cycle}. A cycle $1\dots,m$ is called a \textit{loop} if $m=1$. Note that as the sequence $i,i$ is a loop of $G$ if and only if $e_{ii}\in E$ we may refer to the edge $e_{ii}$ as the loop. 

The main idea behind an isospectral reduction of a graph $G=(V,E,\omega)$ is that we reduce $G$ to a smaller graph on some subset $S\subset V$. To investigate how the structure of a graph is effected by an isospectral reduction we deliberately limit the type of vertex sets over which we can reduce a graph. Such sets, called base sets, are defined as follows.

\begin{definition}\label{def1}\textbf{(Base Set)}
Let $G=(V,E,\omega)$. A nonempty vertex set $S\subset V$ is a \textit{base set} of $G$ if each cycle of $G$, that is not a loop, contains a vertex in $S$.
\end{definition}

A consequence of this definition is that if $S\subset V$ is a base set of the graph $G=(V,E,\omega)$, the subgraph $G|_{\bar{S}}$ has no cycles except possibly for loops.

For $G=(V,E,\omega)$ we let $b(G)$ denote the set of all base sets of the graph $G$. The idea behind the notion of a base set $S\in b(G)$ is the following. Any walk along edges of $G$ that begins at a vertex in $S$ eventually finds its way to another vertex of $S$ in a finite number of steps if we ignore loops. Therefore, a base set allows us to essentially partition a walk on $G$ into finite paths and cycles that begin and end with vertices of $S$. Such paths and cycles are given the following name.

\begin{definition}\textbf{(Branch Product)}
Suppose $G=(V,E,\omega)$ with base set $S=\{1,\dots,m\}$. Let $\mathcal{B}_{ij}(G;S)$ be the set of paths from vertex $i$ to $j$ or cycles if $i=j$ with no interior vertices in $S$. We call a path or cycle $\beta\in\mathcal{B}_{ij}(G;S)$ a \textit{branch} of $G$ with respect to $S$. We let
$$\mathcal{B}_S(G)=\bigcup_{1\leq i,j \leq m} \mathcal{B}_{ij}(G;S)$$
denote the set of all branches of $G$ with respect to $S$.
\end{definition}

If $\beta=1,\dots,m$ is a branch of $G$ with respect to $S$ and $m>2$ define
\begin{equation}\label{eq0.9}
\mathcal{P}_{\omega}(\beta)=\omega(e_{12})\prod_{i=2}^{m-1}\frac{\omega(e_{i,i+1})}{\lambda-\omega(e_{ii})}.
\end{equation}
For $m=1,2$ let $\mathcal{P}_{\omega}(\beta)=\omega(e_{1m})$. We call $\mathcal{P}_{\omega}(\beta)$ the \textit{branch product} of $\beta$.

Notice that the branch product of any $\beta\in\mathcal{B}_S(G)$ is a rational function of $\lambda$. To \emph{isospectrally reduce} a graph over the set $S \in b(G)$, we replace each branch $\mathcal{B}_{ij}(G;S)$ with a single edge $e_{ij}\in\mathcal{E}$. The following definition specifies the weights of these edges.

\begin{definition}\label{IR}\textbf{(Isospectral Graph Reduction)}
Let $G=(V,E,\omega)$ with base set $S=\{1\,\dots,m\}$. Define the edge weights

\begin{equation}\label{eq1.0}
\mu(e_{ij})=\begin{cases}
\displaystyle{\sum_{\beta\in\mathcal{B}_{ij}(G;S)}\mathcal{P}_\omega(\beta)} & \text{if} \ \ \ \mathcal{B}_{ij}(G;S)\neq\emptyset\\
\ \ \ \ \ 0 & \text{otherwise}
            \end{cases} \ \ \ \text{for} \ \ \ 1\leq i,j\leq m.
\end{equation}
The graph $\mathcal{R}_S(G)=(S,\mathcal{E},\mu)$ where $e_{ij}\in \mathcal{E}$ if $\mu(e_{ij})\neq 0$
is the \textit{isospectral reduction} of $G$ over $S$.
\end{definition}

Observe that $\mu(e_{ij})$ in definition \ref{IR} is the weight of the edge $e_{ij}$ in $\mathcal{R}_S(G)$.

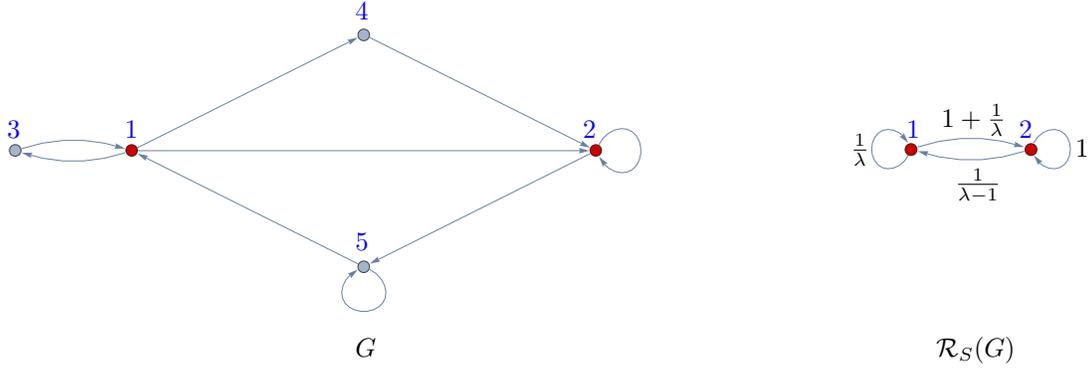
\begin{figure}
  \begin{center}
    \begin{overpic}[scale=.33]{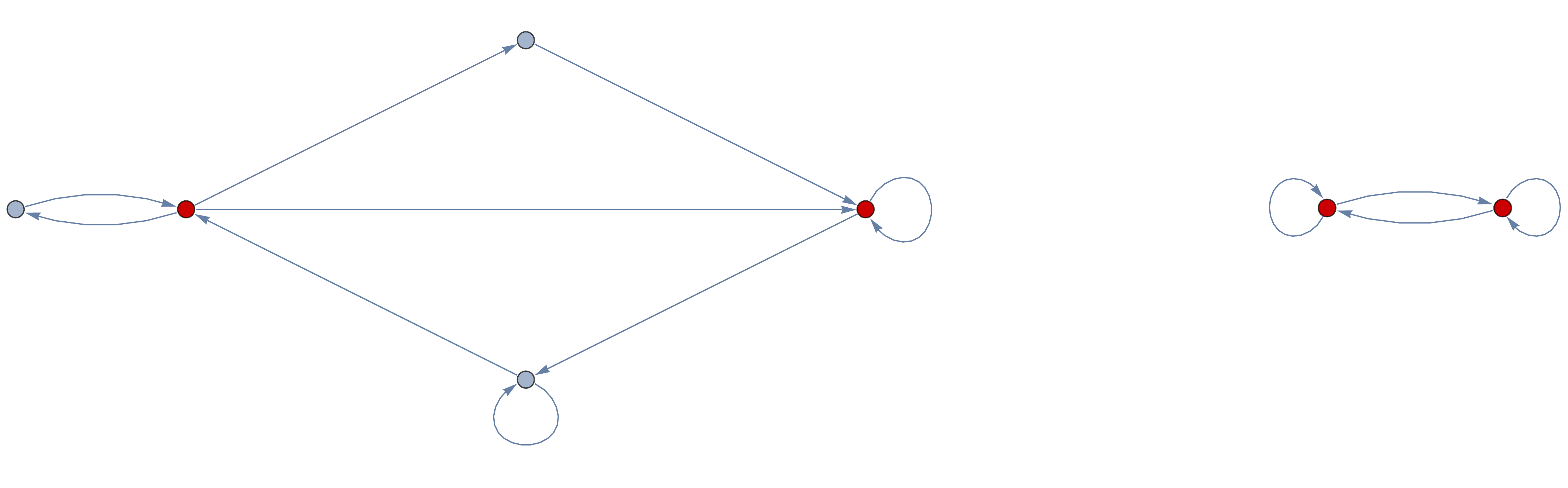}
    \put(32.75,-1.5){$G$}
    \put(11.25,19){\textcolor{blue}{$1$}}
    \put(84.25,19){\textcolor{blue}{$1$}}
    \put(54,19){\textcolor{blue}{$2$}}
    \put(32.75,30){\textcolor{blue}{$4$}}
    \put(32.75,8.5){\textcolor{blue}{$5$}}
    \put(0.25,19){\textcolor{blue}{$3$}}
    \put(94.75,19){\textcolor{blue}{$2$}}
    \put(79,17.2){{$\frac{1}{\lambda}$}}
    \put(87.5,20){$1+\frac{1}{\lambda}$}
    \put(88.75,14){$\frac{1}{\lambda-1}$}
    \put(100,17.2){$1$}
    \put(87,-1.5){$\mathcal{R}_S(G)$}
    \end{overpic}
  \end{center}
  \caption{Reduction of $G$ over $S=\{1,2\}$, shown in red, from example \ref{ex3} where each edge in $G$ has unit weight.}\label{fig2}
\end{figure}

\begin{example}\label{ex3}
Consider the graph $G=(V,E,\omega)$ given in Figure \ref{fig2} (left) where $G$ is an unweighted graph, i.e. each edge of $G$ is given unit weight. Note that the vertex set $S=\{1,2\}\subset V$ is a base set of $G$ since the two nonloop cycles of $G$, namely $1,4,2,5,1$; $1,2,5,1$; and $1,3,1$ each contain a vertex in $S$. In contrast, the vertex set $T=\{2,4,5\}$ is not a base set of $G$ as the (nonloop) cycle $1,3,1$ does not contain a vertex of $T$. Phrased another way, the walk $1,3,1,3,1,\dots$ cannot be partitioned into finite paths and cycles that begin and end with vertices in $T$.

The branches in $\mathcal{B}_S(G)$ are respectively 
\[
\mathcal{B}_{11}(G;S)=\{1,3,1\}, \ \mathcal{B}_{12}(G;S)=\{1,2; \ 1,4,2\}, \ \mathcal{B}_{21}(G;S)=\{2,5,1\}, \ \text{and} \ \mathcal{B}_{22}(G;S)=\{2,2\}.
\]
Using equation (\ref{eq0.9}) the branch product of each branch is given by
$$\mathcal{P}_\omega(1,3,1)=\mathcal{P}_\omega(1,4,2)=\frac{1}{\lambda}, \ \mathcal{P}_\omega(2,5,1)=\frac{1}{\lambda-1} \ \text{and} \ \mathcal{P}_\omega(1,2)=\mathcal{P}_\omega(2,2)=1.$$
Using equation (\ref{eq1.0}) each edge of $\mathcal{R}_S(G)=(S,\mathcal{E},\mu)$ has weight given by
$$\mu(e_{11})=\frac{1}{\lambda}, \ \mu(e_{21})=\frac{1}{\lambda-1}, \ \mu(e_{12})=1+\frac{1}{\lambda}, \ \text{and} \ \mu(e_{22})=1.$$
As each edge weight is nonzero the edge set $\mathcal{E}$ of $\mathcal{R}_S(G)$ is $\mathcal{E}=\{e_{11}$, $e_{12}$, $e_{21}$, $e_{22}\}$. In particular, an edge of $\mathcal{E}$ need not be an edge of $E$. The graph $\mathcal{R}_S(G)$ is shown in Figure \ref{fig2} (right).
\end{example}

Since a graph can only be reduced over specific subsets of its vertex set, i.e. base sets, a natural question is whether a graph can be sequentially reduced to other non-base sets of its vertex set. In \cite{thebook} it has been shown that this is the case where the reduced graph $\mathcal{R}_S(G)$ can be reduced over any base $T\in b(\mathcal{R}_S(G))$ directly using Definition \ref{IR}. The result is the graph $\mathcal{R}_T(\mathcal{R}_S(G))$. In fact the following holds.

\begin{theorem}\label{cor:unique}
\textbf{(Uniqueness of Sequential Reductions)} Let $G=(V,E,\omega)$ be a graph. If $V\supset S_1\supset\cdots\supset S_m$ and $V\supset T_1\supset\cdots\supset T_n \supset S_m$ then the sequence of isoradial reductions
\[
\mathcal{R}_{S_m}(\mathcal{R}_{S_{m-1}}(\cdots\mathcal{R}_{S_1}(G)\cdots))=\mathcal{R}_{S_m}(\mathcal{R}_{T_{n}}(\cdots\mathcal{R}_{T_1}(G)\cdots)).
\]
\end{theorem}

That is, in a sequence of isospectral reductions the result only depends on the final set $S_m$ over which the graph is reduced (see Theorem 2.5 \cite{thebook}). Since the vertex set of a graph with any single vertex removed is a base set of the graph this together with Theorem \ref{cor:unique} has the following consequence. Given a graph $G=(V,E,\omega)$ and subset $S\subset V$ there is a sequence of isospectral reductions that reduces $G$ to a graph on the vertex set $S$. What is more, this reduced graph is unique in that it does not depend on the particular sequence of reductions used to create it. Hence, just as with isospectral matrix reductions it is possible to uniquely reduce a graph over any subset of its vertex set.

This allows us to state the following fundamental relation between isospectral graph and matrix reductions (see Theorem 2.1 in \cite{thebook}).

\begin{theorem}\label{theorem:0}
Let $S$ be a base set of the graph $G=(V,E,\omega)$ with (weighted) adjacency matrix $M$. Then
$$M(\mathcal{R}_{S}(G))=\mathcal{R}_S(M).$$
\end{theorem}

The following result is a corollary to Theorem \ref{theorem:0} and Lemma \ref{lem:symm}, where we allow undirected graphs to possibly have loops.

\begin{theorem}\label{thm:undirect}
Suppose $S$ is a base set of the undirected graph $G=(V,E,\omega)$. Then the reduced graph $\mathcal{R}_S(G)$ is an undirected graph.
\end{theorem}

This theorem will be helpful later on when we prove our main results with latent automorphisms in undirected graphs (as opposed to directed graphs).

\subsubsection{Latent Automorphisms}

We will now describe what we mean by a ``latent automorphism" which generalizes the notion of symmetry in a graph, and which we will use later to characterize cospectral vertices in undirected graphs.

\begin{definition}\label{def:ls}\textbf{(Latent Automorphism)}
We say a graph $G$ has a \emph{latent automorphism} if there exists a subset of vertices which are automorphic in \emph{some} isospectral reduction $\mathcal{R}_S(G)$. If two vertices $a,b\in S$ are automorphic in $\mathcal{R}_S(G)$ we say that $a$ and $b$ are \emph{latently automorphic}. 
\end{definition}

Note here that vertices that are automorphic under the standard definition of \emph{automorphic} are also latently automorphic. The reason is that we can choose the base set $S$ to be all the vertices in the graph, which leaves the graph unchanged by the reduction, i.e. $\mathcal{R}_V(G)=G$.  


\begin{lemma}\label{lem:latent}
Vertices $a$ and $b$ of an undirected graph $G=(V,E,\omega)$ are latently automorphic if and only if $a$ and $b$ are automorphic in $\mathcal R_{\{a,b\}}(G)$
\end{lemma}

\begin{proof}
\noindent($\Leftarrow$) By definition, $a$ and $b$ are latently automorphic if they are are automorphic in  $\mathcal{R}_S(G)$ for some $S\subset V$. Thus this direction is immediate since $\{a,b\}\subset V$.
\smallskip

\noindent($\Rightarrow$) It is easier here to work with the matrices.  Saying that two vertices are automorphic in a graph,  is equivalent to saying that associated adjacency matrix is unchanged after swapping the associated rows and columns.  Suppose $a,b$ are automorphic vertices of $G$, or in other words swapping the associated rows and columns of $M$ leaves M unchanged. For notational simplicity, let $T=\{a,b\}$.  In order to show that $a,b$ are automorphic in $\mathcal{R}_T(G)\in \mathbb{W}^{2\times 2}$, essentially we must show that it is unchanged by interchanging the two rows and the two columns.  In other words we need to show that $$\mathcal{R}_T(M)=\left( \begin{matrix} 0 &1\\1 &0 \end{matrix} \right)\mathcal{R}_T(M)\left(\begin{matrix} 0 &1\\1 &0 \end{matrix}\right).$$
Now we note that because $a,b$ are automorphic in $M$ then $$M_{TT}=\left( \begin{matrix} 0 &1\\1 &0 \end{matrix} \right)M_{TT}\left(\begin{matrix} 0 &1\\1 &0 \end{matrix}\right), \qquad M_{T\bar{T}}=\left( \begin{matrix} 0 &1\\1 &0 \end{matrix} \right) M_{T\bar{T}},\qquad  M_{\bar{T}T}= M_{\bar{T}T}\left( \begin{matrix} 0 &1\\1 &0 \end{matrix} \right).$$ 
Now we simply use the definition of the \emph{Isospectral Matrix Reduction}, to show that 
\begin{align*}
\mathcal{R}_T(M)=& M_{TT}-M_{T\bar{T}}(M_{\bar{T}\bar{T}}-\lambda I)^{-1}M_{\bar{T}T} \\=&\left( \begin{matrix} 0 &1\\1 &0 \end{matrix} \right)M_{TT}\left(\begin{matrix} 0 &1\\1 &0 \end{matrix}\right)-\left( \begin{matrix} 0 &1\\1 &0 \end{matrix} \right) M_{T\bar{T}}(M_{\bar{T}\bar{T}}-\lambda I)^{-1} M_{\bar{T}T}\left( \begin{matrix} 0 &1\\1 &0 \end{matrix} \right) \\=&
\left( \begin{matrix} 0 &1\\1 &0 \end{matrix} \right)\left(M_{TT}-M_{T\bar{T}}(M_{\bar{T}\bar{T}}-\lambda I)^{-1}M_{\bar{T}T}\right) \left(\begin{matrix} 0 &1\\1 &0 \end{matrix}\right)\\=&
\left( \begin{matrix} 0 &1\\1 &0 \end{matrix} \right)\mathcal{R}_T(M)\left(\begin{matrix} 0 &1\\1 &0 \end{matrix}\right).
\end{align*}
Thus $a,b$ are automorphic in $\mathcal{R}_T(G)$.  

By definition, two vertices are latently automorphic if $a,b \text{ are automorphic in } \mathcal{R}_S(G) \text{ for some }S\subset V$.  We can apply the above result to the matrix $\mathcal{R}_S(M)$.  Suppose $a,b$ are automorphic in $\mathcal{R}_S(M)$, then they are also automorphic in $\mathcal{R}_T(\mathcal{R}_S(M))=\mathcal{R}_T(M)$ (this last equality is proven in Theorem 1.3 in \cite{thebook}). Therefore we get the result for the associated graphs for these matrices. 
\end{proof}

An idea that will be discussed in Section \ref{sec:ex} is the notion of \emph{measure of latency} which describes the extent to which the symmetry is ``hidden" within the graph structure, i.e. how small a set do we need to reduce over to see the latent automorphism.  The previous lemma tells us that once a latent automorphism between two vertices appears, the symmetry persists in any reduction over a smaller set.


\section{Main Results}

In this section, we will give the proofs of the main results.  We begin with our characterization of cospectral vertices in terms of isopsectral reductions.  We then will interpret this result in terms of the combinatorics of counting walks in the graph.  Finally, we will give a characterization of strongly cospectral vertices using isospectral reductions.

\subsection{Cospectral Vertices}\label{sec:cosp}

In this section, we prove the main result, which characterizes cospectral vertices in terms of latent automorphisms. It will be useful in our proofs to talk about the isospectral reduction over a single vertex.  We give this the special name of \emph{smash function}.

\begin{definition}\textbf{(Smash Functions)}
The \emph{smash function} of vertex $a\in G$, denoted $S\!F_G(a)$, is the rational function which is the result of performing an isospectral reduction on a graph over the single vertex $a$. That is,
\[S\!F_G(a)=\mathcal{R}_{\{a\}}(G).\]  When $G$ is clear from context, we will suppress it and simply write $S\!F(a)$.
\end{definition}

Our first step is to give the following lemma relating an isospectral reduction over two vertices to the smash functions of those vertices.

\begin{lemma}\label{lem:smash}
For the graph $G=(V,E,\omega)$ with (weighted) adjacency matrix $M=M(G)$ suppose $a,b\in V$. Then $\mathcal R_{\{a,b\}}(M)_{a,a} = \mathcal R_{\{a,b\}}(M)_{b,b}$ if and only if $S\!F_G(a)=S\!F_G(b)$.
\end{lemma}
\begin{proof}

First, suppose  $\mathcal R_{\{a,b\}}(M)_{a,a} = \mathcal R_{\{a,b\}}(M)_{b,b}$. Then $\mathcal{R}_{\{a,b\} }(G)$ must have the form in Figure \ref{fig} and adjacency matrix $$ \left(\begin{matrix}h(\lambda) & f_1(\lambda) \\ f_2(\lambda) & h(\lambda) \end{matrix}\right).$$
Now it is easy for us to calculate $S\!F(a)$ using 
Equation (\ref{IMR}) to be
$$
S\!F(a)=\mathcal{R}_a(M)=M_{11}-M_{12}(M_{22}-\lambda)^{-1}M_{21}=\frac{f_1(\lambda)f_2(\lambda)}{\lambda-h(\lambda)}+h(\lambda)=\mathcal{R}_b(M)=S\!F(b).
$$
We get the same expression for $S\!F(b)$ since $M_{11}=M_{22}$ in this case. Thus, these two smash functions must be equal.
\begin{figure}
\begin{subfigure}[b]{.5\textwidth}
\begin{center}
\begin{overpic}[scale=.2]{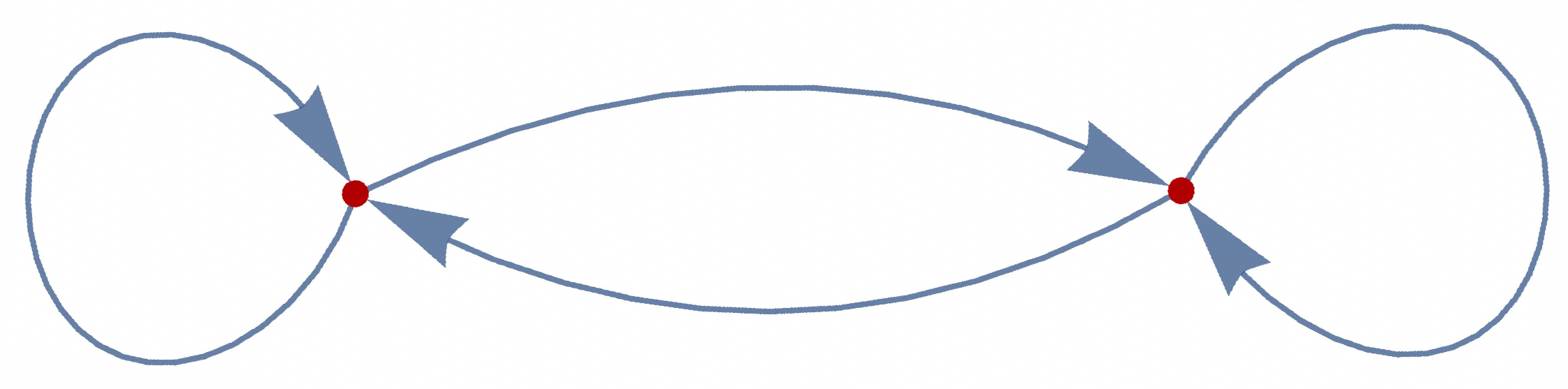}
\put(22,17){$a$}
\put(73,17){$b$}
\put(45,22){$f_1(\lambda)$}
\put(45,0){$f_2(\lambda)$}
\put(-8,11){$h(\lambda)$}
\put(101,11){$h(\lambda)$}
\end{overpic}
\caption{The isospectral reduction over two vertices which are latently automorphic.}\label{fig}
\end{center}
\end{subfigure}~~~~~~
\begin{subfigure}[b]{.5\textwidth}
\begin{center}
\begin{overpic}[scale=.2]{ab1.pdf}
\put(22,17){$a$}
\put(73,17){$b$}
\put(45,22){$f_1(\lambda)$}
\put(45,0){$f_2(\lambda)$}
\put(-8,11){$h(\lambda)$}
\put(101,11){$g(\lambda)$}
\end{overpic}
\caption{The isospectral reduction of an over ANY two vertices. }\label{fig1}
\end{center}
\end{subfigure}
\caption{}\label{fig:red-aut}
\end{figure}



We prove the reverse direction by contradiction. Suppose that  $\mathcal R_{\{a,b\}}(M)_{a,a} \neq \mathcal R_{\{a,b\}}(M)_{b,b}$. When we reduce over $\{a,b\}$, then the graph $\mathcal{R}_{\{a,b\} }(G)$ must have the form of the graph shown in Figure \ref{fig1}, where ${g(\lambda),h(\lambda),f_1(\lambda),f_2(\lambda)\in \mathbb{W}}$, and $h\ne g$.  Now using using 
Equation \eqref{IMR} we calculate the resulting smash functions, suppressing each functions dependence on dependence on $\lambda$ for simplicity. 
$$S\!F(a)=\frac{{{f_1f_2} + h(\lambda - g)}}{{\lambda - g}}$$ $$S\!F(b)=\frac{{{f_1f_2} + g(\lambda - h)}}{{\lambda - h}} $$  
Now suppose by way of contradiction that both vertices do have the same smash function.  Under this supposition, after we set the two expressions equal to each other and simplify we have 
$$\left( {{\lambda^2} - \lambda g - \lambda h + gh - {f_1f_2}} \right)\left( {g - h} \right) = 0.$$  
Because $h \ne g$ then 
$$
{{\lambda^2} - \lambda g - \lambda h + gh - {f_1f_2}}= 0
$$
which implies that  $$g = \frac{{{f_1f_2}}}{{h - \lambda}} + \lambda.$$
As remarked after Definition \ref{def:isored} the rational functions $g$, $h$ and $f_1$, and $f_2$ must have a denominator whose degree is greater than or equal to the degree of their numerator. 


We then calculate 
$$
g\left( \lambda \right) = \frac{f_1f_2}{{h - \lambda}} + \lambda = \frac{{(a_1a_2/b_1b_2)}}{{c/d-\lambda}} + \lambda = \frac{a_1a_2d}{b_1b_2(c-d\lambda)}+\lambda=\frac{p}{q}+\lambda=\frac{p+q\lambda}{q}
$$ 
where $a_1,a_2,b_1,b_2,c,d,p,q$ are polynomials in $\lambda$, $\deg{(a_i)}\leq\deg{(b_i)}$ for $i=1,2$, and $\deg{(c)}\leq\deg{(d)}$. Now if we can show that $\deg{(p)}<\deg{(q)}$ then this would imply $\deg{(p)}\ne \deg{(q+1)}$ and thus 
$$
\deg(p+\lambda q)=\max\{\deg(p),\deg({q}+1)\}>\deg(q).
$$
This contradicts the fact that the degree of the polynomial in the denominator of $g(\lambda)\in\mathbb{W}$ must be greater than or equal to the degree of the polynomial in its numerator.  Thus all that we are required to show is that $\deg{(p)}<\deg{(q)}$, where $p=a_1a_2d$ and $q=b_1b_2(c-d\lambda)$.

Note we know that 
$$\deg{(c)}\leq\deg{(d)}\Rightarrow \deg{(c)}<\deg{(\lambda d)}\Rightarrow\deg{(c-\lambda d)}=\deg{(\lambda d)}>\deg{(d)}.$$  
Since $\deg{(a_i)}\leq\deg{(b_i)}$ for $i=1,2$ and $\deg{(d)}<\deg{(c-\lambda d)}$, we can conclude that 
$\deg{(a_1a_2d)}<\deg{\left(b_1b_2(c-d\lambda)\right)}$, and therefore $\deg{(p)}<\deg{(q)}$. This completes our proof as it follows that $S\!F(a)\neq S\!F(b)$.

\end{proof}

With this lemma, we can now prove our main result for general graphs (both directed and undirected). 

\begin{theorem}\label{thm:cosp}\textbf{(Characterization of Cospectrality via Isospectral Reduction)} 
Let $G=(V,E,\omega)$ with adjacency matrix $M=M(G)$.  Two vertices $a,b\in V$ are cospectral if and only if $\mathcal R_{\{a,b\}}(M)_{a,a}=\mathcal R_{\{a,b\}}(M)_{b,b}$.
\end{theorem}

\begin{proof}
By Lemma \ref{lem:latent} and Lemma \ref{lem:smash}, we need to prove that the vertices $a$ and $b$ are cospectral if and only if $S\!F_G(a)=S\!F_G(b)$.

Using an identity relating determinants and the Schur complement (see \cite{thebook} page 7) we have that
$$
\det(\mathcal{R}_S(M)-\lambda I)=\frac{\det(M-\lambda I)}{\det(M_{\bar{S}\bar{S}}-\lambda I)}
$$
for any $S\subset V$. Thus, if we take $S=\{a\}$ to be the single vertex $a$ we have
\[
\det(\mathcal{R}_{\{a\}}(M)-\lambda I)=\frac{\det(M-\lambda I)}{\det(M_{\bar{a}\bar{a}}-\lambda I)}
\] which gives
\begin{align}\label{eq:smash}
S\!F(a)=\frac{\det(M-\lambda I)}{\det(M_{\bar{a}\bar{a}}-\lambda I)}+\lambda=\frac{p(M,\lambda)}{p(M\backslash a,\lambda)}+\lambda.\end{align}
Now suppose $S\!F(a)=S\!F(b)$ for two vertices $a,b\in G$, then 
$$S\!F(a)=S\!F(b)\iff \frac{p(M,\lambda)}{p(M\backslash a,\lambda)}+\lambda=\frac{p(M,\lambda)}{p(M\backslash b,\lambda)}+\lambda \iff p(M\backslash a,\lambda)=p(M\backslash b,\lambda).$$
Because the characteristic polynomials of $G\backslash a$ and $G\backslash b$ are equal, these subgraphs have the same spectrum. Therefore, the smash functions for two vertices are equal if and only if the vertices are cospectral.

\Hidden{
First, assume that $\RS(M)$ has an automorphism interchanging $u$ and $v$.  Then $\RS(M)(u,u) = \RS(M)(v,v)$.  By Theorem \ref{thm:walks}, this means that
\begin{align}\label{eq:restr_walks}
    |\W_k^*(u,u)| = |\W_k^*(v,v)|
\end{align}
for all $k$.  To prove that $u$ and $v$ are cospectral, by Theorem (FILL REF TO THM) we need to show that
\begin{align}\label{eq:walks}
    |\W_k(u,u)| = |\W_k(v,v)|
\end{align}
for all $k$.  

Any walk from a vertex to itself of length $k$ can be decomposed as a sequence of walks from $u$ to itself that hit $u$ only at the first and last vertex, whose lengths add up to $k$.
Then by (\ref{eq:restr_walks}), we have
\[
|\W_k(u,u)| = \sum_{\ell_1+\cdots+\ell_j=k}\prod_{i=1}^j|\W_{\ell_i}^*(u,u)| = \sum_{\ell_1+\cdots+\ell_j=k}\prod_{i=1}^j|\W_{\ell_i}^*(v,v)| = |\W_k(v,v)|
\]
giving us (\ref{eq:walks}) as required.
}
\end{proof}

This theorem states that when two vertices are cospectral, then the isospectral reduction over those two vertices necessarily has the form in Figure \ref{fig}.  This leads immediately to the following corollary.

\begin{corollary}
If $G=(V,E,\omega)$ is an undirected graph, then the vertices $a,b\in V$ are cospectral if and only if they are latently automorphic.
\end{corollary}
\begin{proof}
This follows directly from Theorem \ref{thm:cosp}, Lemma \ref{lem:symm}, Lemma \ref{lem:latent}, and the observation that a $2\times2$ symmetric matrix has an automorphism if and only if its diagonal entries are equal.
\end{proof}


\subsection{Walk Generating Functions}\label{sec:walk}
In this section, we give a combinatorial interpretation of the isospectral reduction $\RS(M)$ in terms of generating functions of certain walks in the graph $G=G(M)$ associated with $M$.  This will allow us to interpret the results from the previous section combinatorially in terms of walk counts at vertices of $G$.

\begin{definition}\label{def:non-ret_walk}
Given a subset $S\subset V$, we define an \emph{$S$-non-returning} walk in $G$ to be a walk in $G$ that starts and ends in $S$, but all other vertices in the walk (the \emph{internal vertices}) are in $\bar S$.  We will denote by $w^*_k(S)$ the $|S|\times |S|$ matrix whose $(a,b)$ entry is the number of $S$-non-returning walks of length $k$ from $a$ to $b$.

If the graph $G$ has weights on the edges, then this entry is the weighted number of walks.
\end{definition}

\begin{theorem}\label{thm:walks}
Let $M\in\mathbb{C}^{n\times n}$ with corresponding graph $G=(V,E,\omega)$ and let $S\subset V$.  Then 
\[
\RS(M) = \sum_{j=0}^\infty w_{j+1}^*(S) t^j\in \mathbb{W}^{|S|\times |S|}
\]
where $t=1/\lambda$.
\end{theorem}
\begin{proof}
We have
\begin{align*}
    \RS(M) &= M_{SS} - M_{S\bar S}(M_{\bar S\bar S}- \lambda I)^{-1}M_{\bar SS}\\
    &= M_{SS}+\frac1\lambda M_{S\bar S}\left(I-\frac1\lambda M_{\bar S\bar S}\right)^{-1}M_{\bar S S}\\
    &= M_{SS}+t M_{S\bar S}\left(\sum_{k=0}^\infty M_{\bar S\bar S}^kt^k\right)M_{\bar S S}\\
    &= M_{SS}+\sum_{k=0}^\infty M_{S\bar S}M_{\bar S\bar S}^k M_{\bar SS}t^{k+1}\\
    &= M_{SS}+\sum_{j=1}^\infty M_{S\bar S}M_{\bar S\bar S}^{j-1}M_{\bar SS}t^j.
\end{align*}
The result then follows by observing that $M_{\bar S\bar S}^{j-1}$ counts walks of length $j-1$ in $\bar S$, multiplication by $M_{S\bar S}$ counts single steps from $S$ to $\bar S$, and $M_{\bar SS}$ counts single steps from $\bar S$ back to $S$.  Finally, $M_{SS}$ counts walks of length 1 from $u$ to $v$, which trivially satisfy the condition on being in $\bar S$ except at the first and last step.  Thus this can be combined with the sum to give the stated result.
\end{proof}



Recall from Lemma \ref{lem:cosp_char} that one way of characterizing cospectral vertices is that $a$ and $b$ are cospectral if and only if the number of closed walks at $a$ and $b$ of any length are the same.  The previous theorem relates isospectral reductions to non-returning walks at $a$ and $b$.  In what follows, we will connect these two ideas.

We will set some notation.  Recall that $W(M,t)$ denotes the walk generating function for walks given by $M$.  Let us write $W_a(M,t) = W(M,t)_{a,a}$ for the diagonal entry at $a$.  So $W_a(M,t)$ is the generating function for walks beginning and ending at $a$.  That is 
\[
W_a(M,t) = \sum_{\ell=0}^\infty w_\ell(a)t^\ell
\]
where $w_\ell(a)$ denotes the number of (possibly weighted) walks of length $\ell$ from $a$ to itself.  Let us denote by $W_a^*(M,t)$ the generating function for walks beginning and ending at $a$ that do not visit $a$ at any internal step.  That is
\[
W_a^*(M,t) = \sum_{\ell=1}^\infty w_\ell^*(a)t^\ell
\]
where $w_\ell^*(a)$ denotes the number of (possibly weighted) walks of length $\ell$ from $a$ it itself that do not return to $a$ except at the last step (note that we start this sum at $\ell=1$).  Then by Theorem \ref{thm:walks}, we have that
\[
S\!F_G(a)=\mathcal R_{\{a\}}(M) = \lambda W_a^*\left(M,\frac1\lambda\right),
\]
recalling that $t=1/\lambda$.  By Equation (\ref{eq:smash}) in the proof of Theorem \ref{thm:cosp}, we then have
\[
W_a^*\left(M,\frac1\lambda\right) = \frac1\lambda\left(\frac{p(M,\lambda)}{p(M\backslash a,\lambda)}+\lambda\right)
\]
and by Equation (\ref{eq:walk}) of Section \ref{sec:pelim} we have that 
\[
W_a\left(M,\frac1\lambda\right) = -\lambda\frac{p(M\backslash a,\lambda)}{p(M,\lambda)}.
\]
Combining these, we obtain the relationship between these generating functions
\begin{align}
    W_a(M,t) = \frac{1}{1-W_a^*(M,t)}.
\end{align}
Formally manipulating these generating functions, this gives
\[
W_a(M,t) = \sum_{m=0}^\infty W^*_a(M,t)^m = \sum_{m=0}^\infty \left(\sum_{\ell=1}^\infty w_\ell^*(a)t^\ell\right)^m.
\]
From this we obtain the combinatorial identity for any $\ell$
\begin{align}
w_\ell(a) = \sum_{i_1+\cdots+i_k=\ell}\prod_{j=1}^kw_{i_j}^*(a)
\end{align}
where the sum is taken over all tuples partitions $(i_1,...,i_k)$ of non-negative numbers adding up to $\ell$ of any number of parts $k$.  This last identity expresses the number of closed walks at $a$ in terms of the numbers of closed walks at $a$ that do not revisit $a$ until the very last step.  This identity can also be verified in a direct combinatorial way by decomposing a closed walk at $a$ at every return to $a$.  This leads to an alternate way of proving Theorem \ref{thm:cosp}. Thus, another way of thinking of Theorem \ref{thm:cosp} is that the walk counts at $a$ and $b$ are equal for all lengths $\ell$ if and only if the non-returning walk counts at $a$ and $b$ are equal for all $\ell$. 

\subsection{Strongly Cospectral Vertices}\label{sec:strcosp}

In this section, we give a characterization for strongly cospectral vertices (see Definition \ref{def:str-cosp}) in terms of the isospectral reduction.  One way to characterize strongly cospectral vertices is by the following lemma from \cite{godsil2017}.


\begin{lemma}[Corollary 8.4 in \cite{godsil2017}]\label{lem:strcosproots}
Distinct vertices $a,b\in V$ of a graph $G=(V,E,\omega)$ are strongly cospectral if and
only if they are cospectral and all poles of the rational function \[\frac{p(G \backslash\{a, b\}, \lambda)}{p(G, \lambda)}\] are simple.
\end{lemma}


Since $p(G\backslash\{a,b\},\lambda)$ is the characteristic polynomial of the complement of $\{a,b\}$ in $G$, and this comes up in the computation of eigenvalues of isospectral reductions, this allows us to make the connection with the isospectral reduction.

\begin{theorem}\label{thm:str-cosp} \textbf{(Characterization of Strong Cospectrality via Isospectral Reduction)}
Let $G=(V,E,\omega)$ with adjacency matrix $M=M(G)$. Two vertices $a$ and $b$ of $G$ are strongly cospectral if and only if they are cospectral and the eigenvalues of $\RS(M)$ are simple, for $S=\{a,b\}$.
\end{theorem}
\begin{proof}
The eigenvalues of $\RS(M)$ are exactly the roots of the rational function
\[
\frac{p(G,\lambda)}{p(G \backslash\{a, b\},\lambda)}.
\]
These roots are simple if and only if the poles of the reciprocal  \[\frac{p(G \backslash\{a, b\}, \lambda)}{p(G, \lambda)}\] are simple.  Then Lemma \ref{lem:strcosproots} gives the result.
\end{proof}

We remark that for Definition \ref{def:str-cosp} to make sense, the (weighted) adjacency matrix must be symmetric, i.e. the graph in question must be undirected.  However, the isospectral reduction and its eigenvalues can be done for any matrix, symmetric or not.  Thus, Theorem \ref{thm:str-cosp} can be viewed as extending the definition of \emph{strong} cospectrality to directed graphs. 


\section{Applications and Open Questions}\label{sec:ex}

In this section, we will look at some applications of the results that we have proven, and discuss directions for further research.

\subsection{Unpacking an Isospectral Reduction}\label{sec:unpack}
The goal of this section will be to investigate ways of constructing new examples of graphs with cospectral vertices by starting with a graph on two vertices and edges weighted by elements of $\mathbb{W}$ and an automorphism, and then reverse engineering the reduction rules---i.e. to ``unpack" the isospectral reduction---to obtain a graph with a pair of cospectral vertices. 

We will begin with the simplest scheme for unpacking an isospectral reduction.  Start with a reduced graph of the form 
\begin{center}
\begin{overpic}[scale=.2]{ab1.pdf}
\put(22,17){\textcolor{blue}{$1$}}
\put(73,17){\textcolor{blue}{$2$}}
\put(40,22){$f_{12}(\lambda)$}
\put(40,-3){$f_{21}(\lambda)$}
\put(-20,11){$f_{11}(\lambda)$}
\put(101,11){$f_{22}(\lambda)$}
\end{overpic}
\end{center}
whose corresponding matrix is
\begin{equation}\label{eq:last}
\begin{bmatrix}f_{11}(\lambda)&f_{12}(\lambda)\\f_{21}(\lambda)&f_{22}(\lambda) \end{bmatrix}\in\mathbb{W}^{2\times 2}.
\end{equation}
Each $f_{ij}$ is a rational function in $\lambda$ with degree of the numerator less than the degree of the denominator (see the remark after Definition \ref{def:isored}). Then apply the following procedure:\\

\noindent \emph{Step 1}: For each $f_{ij}$, factor the denominator into linear factors (possibly over the complex numbers).\\

\noindent \emph{Step 2}: Now to each of these functions, since the degree of the numerator is at most the degree of the denominator, we can apply partial fraction decomposition and write each function in the form
\[
f_{ij}(\lambda) = \sum_{k=1}^s \frac{a_k}{(\lambda - r_k)^{\ell_k}} \text{  for  } i,j\in\{1,2\}
\]
where the $a_k$ and $r_k$ are (possibly complex) scalars constants, and the $\ell_k$ are non-negative integers.\\

\noindent \emph{Step 3}: Now construct a (complex weighted, directed) graph with vertices $1$ and $2$ as follows.   For each $f_{ij}$, in its decomposition from the previous step, for each $k=1,...,s$, do the following:
\begin{enumerate} \item Attach a directed path $p_{ij}(k)$ from vertex $i$ to vertex $j$ with $\ell_k$ interior vertices.  If $i=j$, this new path starts and ends at the same vertex. If $\ell_k=0$, this simply connects $i$ and $j$ by an edge, which is a loop if $i=j$.
 \item On each new vertex of $p_{ij}(k)$, attach a loop with weight $r_k$.
 \item Choose any edge of the path $p_{ij}(k)$ and weight this with the value $a_k$.  All other new edges receive a weight of 1.
\end{enumerate}

\begin{figure}
  \begin{center}
    \begin{overpic}[scale=.2]{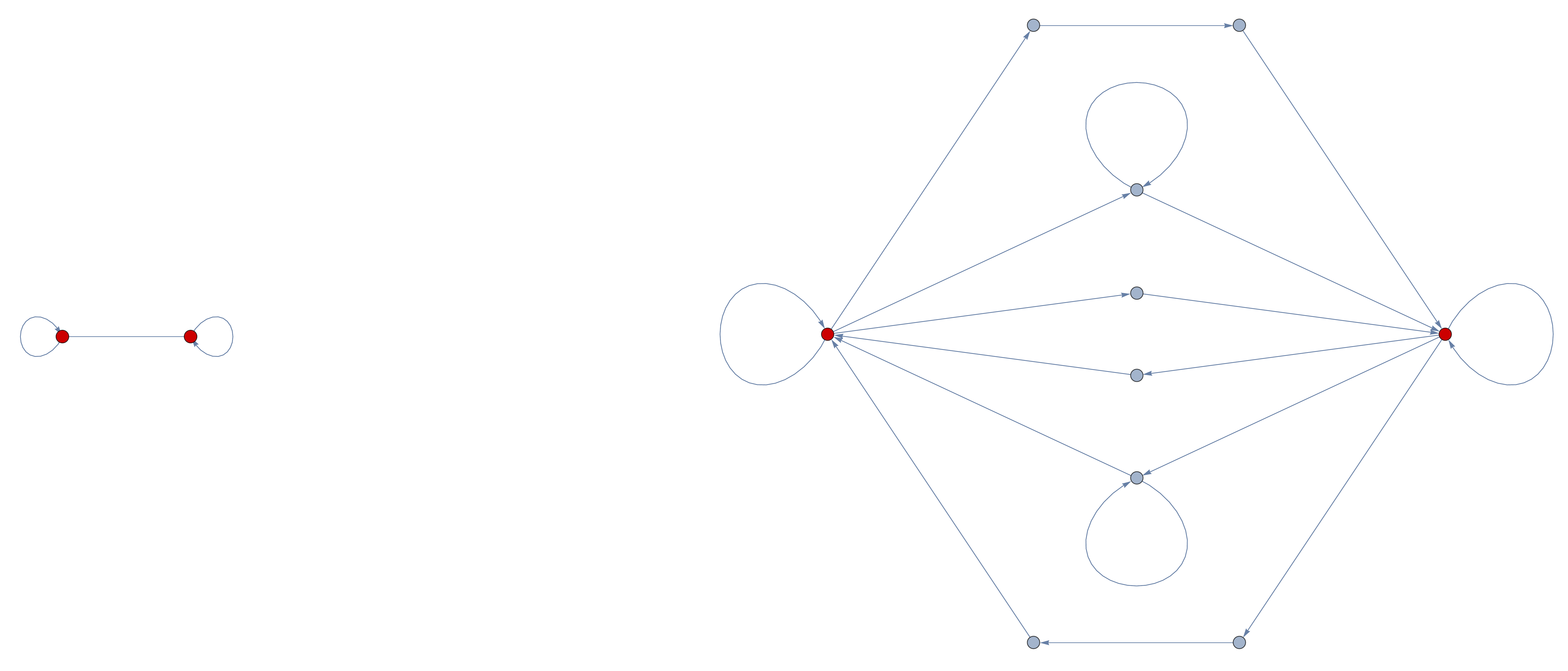}
    \put(71.5,-4){$G$}
    \put(3,10){$\mathcal{R}_{\{1,2\}}(G)$}
    \put(3.25,18){\textcolor{blue}{$1$}}
    \put(52.1,18){\textcolor{blue}{$1$}}
    \put(11.3,18){\textcolor{blue}{$2$}}
    \put(91.5,18){\textcolor{blue}{$2$}}
    \put(5,23.5){{$\frac{1}{\lambda^3-\lambda^2}$}}
    \put(-1,20.5){$1$}
    \put(15.5,20.5){$1$}
    \put(43.5,20.5){$1$}
    \put(100,20.5){$1$}
    
    \put(55,32){$-1$}
    \put(72,42.5){$1$}
    \put(86.5,32){$1$}
    
    \put(62,27.5){$1$}
    \put(72,38.25){$1$}
    \put(81.5,27.5){$1$}
    
    \put(63,24){$-1$}
    \put(78.5,24){$1$}
    
    \put(55,9){$-1$}
    \put(72,-1){$1$}
    \put(86.5,9){$1$}
    
    \put(62,13.5){$1$}
    \put(72,3){$1$}
    \put(81.5,13.5){$1$}
    
    \put(63,17.25){$-1$}
    \put(78.5,17.25){$1$}
    
    \end{overpic}
  \end{center}
  \caption{The reduction $\mathcal{R}_{\{1,2\}}(G)$ is unpacked using the rules described in Section \ref{sec:unpack}. Since the loops of vertices 1 and 2 have the same weights in $\mathcal{R}_{\{1,2\}}(G)$ then these vertices are cospectral in $G$ although they are not automorphic.}\label{fig5}
\end{figure}
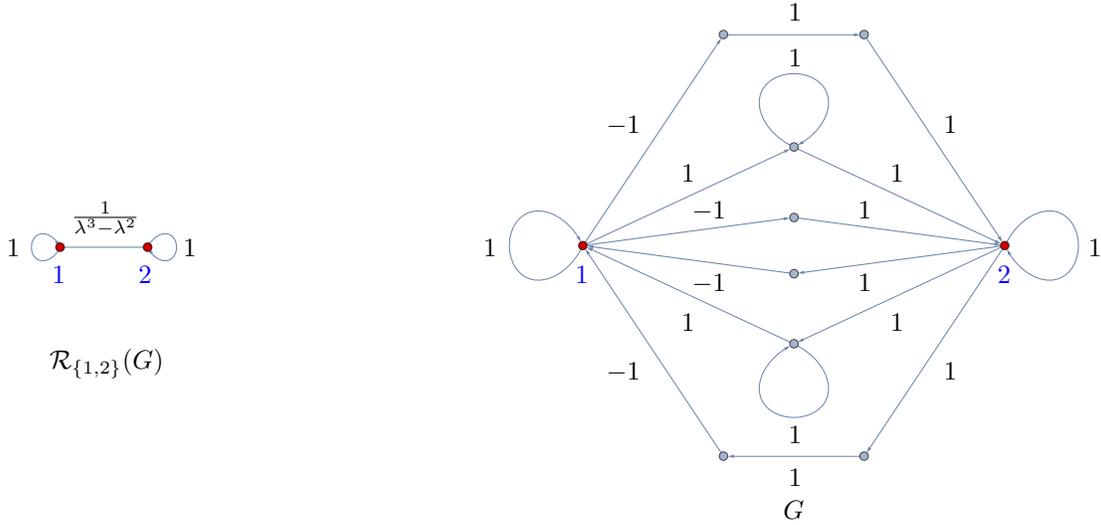

\noindent Then by the isospectral graph reduction rules (see Equations (\ref{eq0.9}) and (\ref{eq1.0})), reducing this graph over vertices $1$ and $2$ will give us the graph we started with, namely the graph with (weighted) adjacency matrix in given in Equation \eqref{eq:last}.  If $f_{11} = f_{22}$, then Theorem \ref{thm:cosp} guarantees that vertices $1$ and $2$ will be cospectral.  

This construction has the advantage of being relatively straightforward, but has the drawback that the weights, in general, can be complex numbers. The graph is also by construction a directed graph with a corresponding non-symmetric adjacency matrix.  But this method does give a straightforward way of producing a complex-weighted directed graph with a pair of cospectral vertices.  

To give a concrete example of this process consider the undirected graph $\mathcal{R}_{\{1,2\}}$ in Figure \ref{fig5} with weighted adjacency matrx
\[
W=\begin{bmatrix}1&\frac{1}{\lambda^3-\lambda^2}\\\frac{1}{\lambda^3-\lambda^2}&1 \end{bmatrix}\in\mathbb{W}^{2\times 2}.
\]
Following Steps 1 and 2 of our unpacking process the partial fraction decomposition of the entries of $W_{ij}=f_{ij}$ are 
\[
f_{11}=f_{22}=1 \text{ and } f_{12}=f_{21}=\frac{1}{\lambda-1}+\frac{-1}{\lambda^2}+\frac{-1}{\lambda}.  
\]
Hence, the loops of vertices 1 and 2 remain unchanged but the edge from vertex 1 to vertex 2 is replaced with three paths $p_{12}(1)$, $p_{12}(2)$, and $p_{12}(3)$. Path $p_{12}(1)$ has length 2 where the single interior vertex has a loop and all edges have unit weight. Path $p_{12}(2)$ has length 3, no loops, and we choose the first edge to have weight $-1$. Path $p_{12}(3)$ has length 2, no loop(s), and we choose the first edge to have weight $-1$.

We similarly construct the paths $p_{21}(1)$, $p_{21}(2)$, and $p_{21}(3)$ with the exception that we choose the last edges in $p_{21}(2)$ and $p_{21}(3)$ to have weight $-1$. This results in the graph $G$, shown right, in which vertices 1 and 2 are not automorphic but are cospectral since the loops of these vertices have the same weights in the reduced graph $\mathcal{R}_{\{1,2\}}$.

The question of unpacking the isospectral reduction to produce an unweighted, undirected graph appears to be much more subtle.  Consider, for instance, the $\W$-weighted graph in Figure \ref{fig:weighted_square}. 
\tikzstyle{every node}=[circle,draw,fill=white,inner sep=1pt]
\begin{figure}
\begin{subfigure}[b]{.5\textwidth}
    \centering
    \begin{tikzpicture}
     \draw (90:1)node(a){} (210:1)node(b){} (330:1)node(c){};
     \draw (a) to node[draw=none]{$\frac1\lambda$} (b);\draw (b) to node[draw=none]{$\frac1\lambda$} (c); \draw (c) to node[draw=none]{$\frac1\lambda$} (a);
     \draw (a) edge[loop above] node[draw=none]{$\frac2\lambda$} (a);
     \draw (b) edge[loop left] node[draw=none]{$\frac2\lambda$} (b);
     \draw (c) edge[loop right] node[draw=none]{$\frac2\lambda$} (c);
    \end{tikzpicture}
    \caption{A reduction $\RS(G)$ with symmetry}
    \label{fig:weighted_square}
\end{subfigure}
\begin{subfigure}[b]{.5\textwidth}
  \begin{center}
    \begin{tikzpicture}
    \draw \foreach \x in {90,210,330}
    {
      (\x:0)node{}--(\x:1)node[fill=black]{}--(\x:2)node{}
    }
    ;
    \end{tikzpicture}~~~~~~~~~~~~~
    \begin{tikzpicture}
    \draw \foreach \x in {90,210,330}
    {
      (\x:1.5)node[fill=black]{}--(\x+60:1.5)node{}--(\x+120:1.5)
    }
    ;
    \end{tikzpicture}
    \end{center}
    \caption{Two graphs whose isospectral reduction is the same.} 
    \label{fig:non-unique}
\end{subfigure}
\caption{}
\end{figure}
Examining the rational functions along each edge and loop, with the goal of unpacking this to a simple unweighted graph $G$, it becomes apparent that each of the three vertices that get reduced down to the vertices of this graph must be adjacent to two other vertices in $G$ because of the weights on the loops. In addition,  each pair of vertices that are adjacent here must be connected in $G$ by a path with one vertex in it.  Now, we can choose for each of these paths to be different, or they can all include the same vertex, and then have each of the three original vertices connected to its own separate vertex.  Thus we can see that for the two graphs in Figure \ref{fig:non-unique}, the reduction over set of vertices colored black both yield the graph in Figure \ref{fig:weighted_square}. 
This shows us that the ``unpacking" of an isospectral reduction will not, in general, be unique.

As a further example, consider a modification of this example in Figure \ref{fig:unpack}. If we wish to unpack this graph, it is similar to the previous example, but now we have four vertices, but between two of them, there should not be a path in the unpacked graph $G$ not using the original vertices.  This is accomplished by the graph in Figure \ref{fig:new-cosp}---the reduction over the black vertices yields the reduced graph in Figure \ref{fig:unpack}.
Because of the automorphism in the reduced graph, the vertices labelled $a$ and $b$ are cospectral.  This is an example of a graph with a pair of cospectral vertices that does not have an automorphism between them (the other two black vertices are also cospectral, and are automorphic). In fact, vertices $a$ and $b$ are strongly cospectral.  This example does not seem to fall into any category of standard constructions of cospectral vertices or strongly cospectral vertices (see our discussion of constructing cospectral pairs in Section \ref{subsec:prelim_cospectral}).

\begin{figure}
\begin{subfigure}[b]{.5\textwidth}
    \centering
    \begin{tikzpicture}
     \draw (0,-1)node(a){} (1,0)node(b){} (0,1)node(c){} (-1,0)node(d){};
     \draw (a) to node[draw=none]{$\frac1\lambda$} (b); \draw (b) to node[draw=none]{$\frac1\lambda$}(c); \draw (c) to node[draw=none]{$\frac1\lambda$} (d); \draw(d) to node[draw=none]{$\frac1\lambda$} (a); \draw (a) to node[draw=none]{$\frac1\lambda$} (c);
     \draw (a) edge[loop below] node[draw=none]{$\frac2\lambda$} (a);
     \draw (b) edge[loop right] node[draw=none]{$\frac2\lambda$} (b);
     \draw (c) edge[loop above] node[draw=none]{$\frac2\lambda$} (c);
     \draw (d) edge[loop left] node[draw=none]{$\frac2\lambda$} (c);
    \end{tikzpicture}
    \caption{Another reduction with a symmetry}
    \label{fig:unpack}
    \end{subfigure}
    \begin{subfigure}[b]{.5\textwidth}
       \centering
    \begin{tikzpicture}
    \draw \foreach \x in {60,180,300}
    {
    (\x:1)node[fill=black]{}--(\x+60:1)node{}--(\x+120:1)
    }
    (1,0)--(2,0)node[fill=black]{}--(3,0)node{};
    \draw (-1.1,0)node[draw=none,left]{$a$} (2,.1)node[draw=none,above]{$b$};
    \end{tikzpicture}
    \caption{A graph whose reduction is the graph in Figure \ref{fig:unpack}.} 
    \label{fig:new-cosp}
    \end{subfigure}
    \caption{}
\end{figure}

An open question is to determine conditions on the rational function that will guarantee that a reduced matrix comes from a simple graph, and to determine a set of ``unpacking rules" which will construct a simple unweighted graph from the reduced.  Various relaxations of this question that would also be interesting are, if there are conditions and rules that would guarantee unpacking to a weighted undirected graph, and conditions and rules that would guarantee unpacking to an unweighted directed graph.  Even partial answers to any of these questions we feel could have deep implications to graph inverse eigenvalue problems (see for instance \cite{hogben2005spectral}).

\subsection{Measure of Latency}
In \cite{smith2019hidden}, in connection with hidden symmetries and isospectral reduction, the notion of \emph{measure of latency} is introduced.  This concept gives a measure of how ``hidden" the symmetry is when there is a latent automorphism.  With our characterization of cospectrality, every pair of cospectral vertices in an undirected graph are latently automorphic, so this measure of latency can be applied to the theory of cospectral vertices to measure how far a pair of cospectral vertices is from being automorphic.  

\begin{definition}[Measure of Latency]
Let $G=(V,E,\omega)$ be a graph with $n$ vertices and let $S$ be a subset of its vertices which are latently symmetric. The \emph{measure of latency} $\M$ is defined as \[\M(S)=\frac{n-|T|}{n-|S|}\] where $T\subseteq V$ is a maximal set of vertices such that the vertices $S$ are symmetric in $\mathcal R_T(G)$.
\end{definition}

We will consider $S=\{a,b\}$ where $a$ and $b$ are cospectral vertices of $G=(V,E,\omega)$, and the measure of latency involves the largest $T$ containing $a$ and $b$ such that $\mathcal R_T(G)$ has an automorphism permuting $a$ and $b$.  Of course, if there is already an automorphism of $G$ that maps $a$ to $b$, then we take $T$ to be the entire vertex set $V$, and $\M=0$.  If we must reduce all the way to $S=\{a,b\}$ before an automorphism appears, then $\M =1$.  We will compute the measure of latency for the graphs with cospectral vertices mentioned in this paper.

For the graphs in Figures \ref{fig:isom_no_aut} and \ref{fig:non-aut_cosp}, if we work through the possible intermediate reductions over all $T\supset S=\{a,b\}$, we see that only for $T=S$ is there an automorphism between $a$ and $b$. Thus, for each of these graphs, we have $\M = 1$.

\Hidden{ 
\begin{center}
\begin{tikzpicture}
\draw (0,0)node[fill=black]{}--(-1,.6)node[fill=black]{}--(-1,-.6)node{}--(0,0)--(1,1)node{}--(2.1,.6)node{}--(2.1,-.6)node{}--(1,-1)node{}--(0,0)
(1,1)--(1,2)node{}
(2.1,.6)node{}--(3,0)node[fill=black]{}--(2.1,-.6)node[fill=black]{}
(-1.1,.6)node[draw=none, left]{\small $a$}
(3.1,0)node[draw=none, right]{\small $b$};
\end{tikzpicture}~~~~~
\begin{tikzpicture}
\draw (0,0)node[fill=black](a){} (2,0)node[fill=black](b){} (4,0)node[fill=black](c){} (6,0)node[fill=black](d){};
\draw (a) to node[draw=none,above]{$1+\frac1\lambda$} (b); \draw (b) to node[draw=none,above]{$\frac1\lambda+\frac{1}{\lambda^2-2}$} (c); \draw (c) to node[draw=none,above]{$1+\frac1\lambda$} (d);
\draw (a) edge[loop left] node[draw=none]{$\frac1\lambda$} (a); 
\draw (b) edge[loop below] node[draw=none]{$\frac2\lambda$} (b); 
\draw (c) edge[loop below] node[draw=none]{$\frac2\lambda$} (c); 
\draw (d) edge[loop right] node[draw=none]{$\frac1\lambda$} (d); 
\end{tikzpicture}
\end{center}
}
For the graph in Figure \ref{fig:new-cosp}, we already saw that a reduction down to four vertices has an automorphism (Figure \ref{fig:unpack}). Thus,
\[
\M = \frac{8-4}{8-2}=\frac23.
\]

An interesting avenue for future investigation is to study the combinatorial properties of the measure of latency.  From these examples, it appears that a smaller measure of latency is associated with the presence of other automorphisms in the graph. It would be interesting to determine if this holds in general.  It would also be of interest to study the combinatorial meaning of larger subsets being latently automorphic, and if there is any characterization similar to Theorem \ref{thm:cosp} for subsets of size larger than 2.  Finally, it would be interesting to see if there is any relationship between this measure of latency and strong (versus non-strong) cospectrality. 


\bibliographystyle{plain}
\bibliography{bibfile.bib}
\end{document}